# A New Efficient Methodology for AC Transmission Network Expansion Planning in The Presence of Uncertainties

Soumya Das, *Student Member, IEEE*, Ashu Verma, *Sr. Member, IEEE* and P. R. Bijwe, *Sr. Member, IEEE*

*Abstract*—Consideration of generation, load and network uncertainties in modern transmission network expansion planning (TNEP) is gaining interest due to large-scale integration of renewable energy sources with the existing grid. However, it is a formidable task when iterative AC formulation is used. Computational burden for solving the usual ACTNEP with these uncertainties is such that, it is almost impossible to obtain a solution even for a medium-sized system within a viable time frame. In this work, a two-stage solution methodology is proposed to obtain quick, good-quality, sub-optimal solutions with reasonable computational burden. Probabilistic formulation is used to account for the different uncertainties. Probabilistic TNEP is solved by 2m+1-point estimate method along with a modified artificial bee colony (MABC) algorithm, for Garver 6 bus and IEEE 24 bus systems. In both the systems, rated wind generation is considered to be more than one-tenth of the total generation capacity. When compared with the conventional single stage and existing solution methods, the proposed methodology is able to obtain almost identical solutions with extremely low computational burdens. Therefore, the proposed method provides a tool for efficient solution of future probabilistic ACTNEP problems with greater level of complexity.

NOMENCLATURE

The primary notations used in this paper are given below. Other symbols used in the paper are defined as per requirement.

A. *Parameters, Variables and Vectors related to TNEP*

| | |
|---|---|
| $v_{cr}$ | Final line investment cost for crisp/deterministic TNEP |
| $v_{pr}$ | Final line investment cost for probabilistic TNEP |
| $l$ | Power corridor between two buses |
| $C_l$ | Cost of installation of a line in the $l^{th}$ power corridor |
| $n_l$ | Number of additional lines in the $l^{th}$ power corridor |
| $\Omega$ | Set of all power corridors |
| $k$ | Contingency state, $k = 0$ denotes base case |
| $P_{Dem}$ | Vector of real power demand |
| $Q_{Dem}$ | Vector of reactive power demand |
| $q_{reac}$ | Vector of capacity of additional reactive sources |
| $P_{inj}^k$ | Vector of real power injections at $k^{th}$ contingency |
| $P_{Gen}^k$ | Vector of real power generations at $k^{th}$ contingency |
| $Q_{inj}^k$ | Vector of reactive power injections at $k^{th}$ contingency |
| $Q_{Gen}^k$ | Vector of reactive power generations at $k^{th}$ contingency |
| $V^k$ | Vector of bus voltage magnitudes at $k^{th}$ contingency |
| $\theta^k$ | Vector of bus voltage angles at $k^{th}$ contingency |
| $n^k$ | Vector of system circuits at $k^{th}$ contingency |
| $S_l^{kfr}$ | Apparent power flow at the sending end of each line of $l^{th}$ corridor at $k^{th}$ contingency |
| $S_l^{kto}$ | Apparent power flow at the receiving end of each line of $l^{th}$ corridor at $k^{th}$ contingency |
| $S_{lmax}$ | Maximum power limit of a line in $l^{th}$ power corridor |
| $\bar{n}_l$ | Allowable limit of new lines in $l^{th}$ power corridor |
| $\omega_K$ | Total number of contingencies |
| MABC | Modified Artificial Bee Colony |
| $v_{Aug}$ | Augmented objective function |
| $F_{pen}$ | Penalty function |
| $F_{Eqpen}$ | Penalty function for equality constraints |
| $F_{Ieqpen}$ | Penalty function for inequality constraints |
| $\gamma_{Eq}$ | Weighting factor for equality constraints |
| $\gamma_{Ir}$ | Weighting factor for inequality constraints |
| $f_x$ | Probability density function (pdf) of a random variable $x$ |
| $\alpha$ | Shape parameter of Weibull function |
| $\beta$ | Scale parameter of Weibull function |
| $U$ | Wind speed as a random variable |
| $L$ | Load at a bus as a random variable |
| $P_{WTG}$ | Wind turbine generator (WTG) output |
| $P_{WTG_{rt}}$ | Rated WTG output |
| $U_{ci}$ | Cut-in wind speed for WTG |
| $U_{rt}$ | Rated wind speed for WTG |
| $U_{co}$ | Cut-out wind speed for WTG |
| $\rho_{FOR}$ | Forced outage rate (FOR) of a line |
| $E(-)$ | Expected value of a random variable |
| $\mu_x$ | Expected value of a random variable $x$ |
| $\sigma_x$ | Standard deviation of a random variable $x$ |

B. *Parameters of the Metaheuristic Algorithm*

| | |
|---|---|
| $cs_N$ | Population of generated solutions |
| $\psi$ | Number of neighbours |
| $iter$ | Maximum number of iterations per trial |
| $lim$ | Parameter for generation of scout bees |
| $w_g$ | Weightage to control the effect of best solution on bee movement |
| $tp$ | Approximate time required for a single trial |

## I. INTRODUCTION

GLOBAL warming is a major phenomenon in the present world. In accordance with several environmental conventions [1] laid down over the past couple of decades and most significantly, with the Kyoto protocol [2] and the Paris agreement [3] now in force, major economies are introducing policies for more renewable power generation over conventional power plants. This has led to increased integration of renewable energy sources (RES) into the existing grid. However, due to their inherent intermittency and uncertainty, modern grid operation and planning are facing substantial challenges.

Network planners in the current scenario must consider these uncertainties along with the uncertainties of load and network availability while arriving at a feasible plan looking forward into future grid expansions. Transmission network expansion planning (TNEP) in itself is an extremely difficult task to accomplish due to its mixed-integer, non-linear and non-convex characteristic. Consideration of these uncertainties increases its complexity manifold. Therefore, mostly, when such uncertainties are incorporated, planning is performed with simplistic DC formulation [4]-[15].

Also, robust solution methodologies for DCTNEP are formulated in [14] and [15] for handling uncertainties related to load and generation. Robust formulation does not require probabilistic information of the uncertain variables, only the extreme limits of the variation is enough for its application. This is particularly useful where probabilistic distribution of a variable if difficult to obtain. However, in TNEP, non-consideration of probability of occurrence of a particular state of an uncertain variable may lead to over estimation of the system resources, thereby increasing the planning cost.

DCTNEP when applied to an actual AC system is extremely likely to cause several operational problems due to its simplistic and approximate model with neglected reactive power flows and respective losses. As a result, modern trend is to solve ACTNEP. However, full, non-linear ACTNEP is computationally extremely demanding due to the requirement of repeated iterative AC power flow (ACPF) solutions. In addition, solution of such non-linear ACTNEP in presence of uncertainties is computationally so intensive that obtaining a solution for even a moderately sized system becomes impractical. Therefore, many researchers have considered linearized AC models for solving the same. Authors in [16] present a mixed integer linear programming (MILP) based methodology to solve probabilistic ACTNEP (PrACTNEP). The base case solution obtained is then checked for N-1 contingencies and suitable reinforcements are added. Although very efficient computationally, this method of planning almost always produces sub-optimal results. Similarly, MILP is used to solve a stepwise long-term investment plan for a pseudo-dynamic PrACTNEP in [17]. Load and generation uncertainties are considered without any component outages.

Even when a non-linear AC modelling is used, several assumptions are considered to reduce the computational burden. Likewise, a new variation of binary differential evolution (BDE) algorithm is proposed for long-term probabilistic AC/DC transmission expansion planning in [18]. However, network component outage probabilities are not considered here. In [19], authors propose a concept that economic adjustment of TNEP should be done in response to extreme events. PrACTNEP is solved through scenario generation by Monte-Carlo simulations (MCS). To reduce the computational burden, only a predetermined set of contingencies are checked in place of deterministic N-1 criterion. Also, computational burden is reduced by performing selective N-1 security analysis for an expansion plan at only peak load condition. Probabilistic analysis is conducted only if there is no network violation in the previous stage. However, as a network is prone to violations due to contingencies even at other loading scenarios [20], the methodology in [19] may produce violations in a real system. Authors in [21] present a reliability based probabilistic approach for ACTNEP. Based on a couple of cost indices, the solution methodology involves successive reinforcement of an existing network in response to all N-1 and selected N-2 contingencies. Such a method of planning is extremely likely to provide much higher investment costs than those obtained through full mixed-integer planning methodologies. MCS based non-dominated sorting DE algorithm is used for PrACTNEP in [22]. Uncertainty of loads, demand response programs (DRPs) and distributed generation (DG) are considered, whereas, deterministic N-1 security criterion is not used in this study.

From the existing literature review, it can be observed that, probabilistic ACTNEP is not a very matured area of research. Most of the methods of solving the same consider one or more simplifications to reduce the computational burden. Further, as per the best of knowledge of the authors, combined solution of non-linear PrACTNEP with deterministic N-1 security analysis has not been attempted in the past. Although such a planning solution is extremely useful in a real power system with several uncertainties, it is avoided mainly due to its extreme computational burden.

Therefore, to address this issue in current planning studies, the authors in this study present an efficient, two-stage, PrACTNEP solution methodology based on $2m + 1$ point estimate method [23], [24] and a modified artificial bee colony algorithm [25]. The primary contributions of this paper can be summarized as: 1) development of an intelligent, two-stage solution methodology for solving PrACTNEP and, 2) its application to solve the same for two test systems, demonstrating its ability to obtain good-quality, sub-optimal solutions with drastically low computational burdens. Remaining sections of this paper are organized as follows: Section II discusses the mathematical modelling, Section III provides a brief overview of point estimate method, while, details of the proposed solution methodology are in Section IV. Results and discussion are provided in Section V, followed by the conclusion of this paper in Section VI.

## II. MATHEMATICAL MODELLING

The mathematical model for probabilistic AC TNEP is based on conventional static, deterministic/crisp AC TNEP, which is discussed in brief as follows:

### A. Conventional static, deterministic ACTNEP:

Modelling of static, deterministic ACTNEP is pretty straight forward as it does not consider any uncertainties [27]:

Minimize:
$$v_{cr} = \sum_{l \in \Omega} C_l \times n_l \quad (1)$$

Subject to:

*1) Network constraints:*

$\forall k,$
$$P_{inj}^k - P_{Gen}^k + P_{Dem} = 0$$
$$Q_{inj}^k - Q_{Gen}^k + Q_{Dem} - q_{reac} = 0$$
$$P_{Genmin} \leq P_{Gen}^k \leq P_{Genmax}$$
$$Q_{Genmin} \leq Q_{Gen}^k \leq Q_{Genmax}$$
$$V_{min} \leq V^k \leq V_{max} \quad (2)$$

Power flows at both ends of individual lines in a corridor, $l$ should be within their respective limits:

$\forall l \in \Omega \;;\; l \neq k,$
$$(n_l^0 + n_l)S_l^{kfr} \leq (n_l^0 + n_l)S_{lmax} \quad (3)$$
$$(n_l^0 + n_l)S_l^{kto} \leq (n_l^0 + n_l)S_{lmax} \quad (4)$$

for $l = k, k \neq 0,$
$$(n_l^0 + n_l - 1)S_l^{kfr} \leq (n_l^0 + n_l - 1)S_{lmax} \quad (5)$$
$$(n_l^0 + n_l - 1)S_l^{kto} \leq (n_l^0 + n_l - 1)S_{lmax} \quad (6)$$

Where, power flow at the sending end is,
$$S_l^{kfr} = \sqrt{\left(P_l^{kfr}\right)^2 + \left(Q_l^{kfr}\right)^2} \quad (7)$$

Power flow at the receiving end is provided as,
$$S_l^{kto} = \sqrt{(P_l^{kto})^2 + (Q_l^{kto})^2} \quad (8)$$

$P_l^{kfr}$, $P_l^{kto}$, $Q_l^{kfr}$ and $Q_l^{kto}$ are the respective real and reactive power flows through individual lines of $l^{th}$ power corridor [28].

*2) Physical constraints:*

There is a maximum limit of the number of new lines that can be constructed in a power corridor. This condition is taken care of by the following constraint:
$$\forall l \in \Omega \;;\; \forall k, \quad 0 \leq n_l \leq \bar{n}_l \quad (9)$$

Here, $n_l \geq 0$ and integer for $\forall l \in \Omega$ and $l \neq k$; for $l = k$, $k \neq 0$, $(n_l - 1) \geq 0$ and integer. $k = 0, 1, \dots \omega_K$. Subscripts $min$ and $max$ are used in the notations of the variables to denote their respective minimum and maximum values.

Equation (1) is the objective function, which is the total line investment cost. The usual constraints to be satisfied for an AC OPF are provided by (2). These are the network real and reactive power equality constraints; power generator production limits; and limitations on the size of additional reactive sources and system voltage magnitudes.

Reactive power planning (RPP) is an integral part of ACTNEP. However, in this work, to reduce the computational burden of solving TNEP, it is assumed that proper RPP is already available. Therefore, $q_{reac}$ is fixed for all network contingencies, and hence, their installation costs are not included in the objective function. Our mathematical modelling also considers that in a power corridor, lines having exactly similar characteristic are installed. Whenever a new line having different characteristic from the existing lines is installed, it is considered to be representing a separate sub-corridor within the same set of system buses. Thus, '$l$' represents such sub-corridors in the entire network. This way of modelling is advantageous in the way that it can track outages of all different types of lines present in a system.

### B. Solution methodology of deterministic ACTNEP by MABC:

As with all metaheuristic algorithms, MABC is only able to solve unconstrained optimization problems. Therefore, to effectively solve ACTNEP by MABC, all the constraints are included within the original objective function as penalty functions. Purpose of such penalty function approach is to add a high value to the original objective function only when there is a constraint violation, otherwise, the value addition is zero. Therefore, when all the constraints are satisfied, the augmented objective function becomes equal to the original. Thus,
$$v_{Aug} = v_{cr} + F_{pen} \quad (10)$$

Where, $\quad F_{pen} = \gamma_{Eq} \times F_{Eqpen} + F_{Ieqpen} \quad (11)$

$F_{Eqpen}$ denotes penalty function for the network equality constraints, with $\gamma_{Eq}$ as the weighting factor. Satisfaction of the network equality constraints set $F_{Eqpen} = 0$.

$F_{Ieqpen}$ is the penalty function for all other inequality constraints that need to be satisfied. It can be defined as:
$$F_{Ieqpen} = \sum_{r=1}^{o} \tau_r \quad (12)$$

Where,
$$\tau_r = \begin{cases} \gamma_{Ir}(|\chi_{rmin}| - |\chi_r^k|)^2 & |\chi_r^k| < |\chi_{rmin}| \\ 0 & |\chi_{rmin}| \leq |\chi_r^k| \leq |\chi_{rmax}| \\ \gamma_{Ir}(|\chi_r^k| - |\chi_{rmax}|)^2 & |\chi_r^k| > |\chi_{rmax}| \end{cases} \quad (13)$$

$\gamma_{Ir}$ denote the weighting factor for $r$-th inequality constraint and, $\chi_{rmin}$ and $\chi_{rmax}$ are their corresponding minimum and maximum limits. Therefore, $F_{Ieqpen} = 0$ only when all the inequality constraints are satisfied. After the evaluation of $v_{Aug}$, fitness value is obtained as:
$$fitness = \frac{1}{v_{Aug}} \quad (14)$$

This fitness value is used by MABC to quantify the quality of a solution obtained and search for global minimum.

### C. Probabilistic ACTNEP:

When uncertainties in a system are considered probabilistically to solve ACTNEP, a probabilistic formulation of (1) is required to be solved. Therefore, for solution by MABC, as with any other metaheuristic, a probabilistic realization of (10) is required to be solved. This is obtained by replacing $v_{cr}$ with line cost for probabilistic TNEP, $v_{pr}$ and crisp penalty function $F_{pen}$ in (10) by its expected value:
$$v_{Aug} = v_{pr} + E(F_{pen}) \quad (15)$$

$v_{pr}$ is obtained in the same way as (1). Expected value of $F_{pen}$ can be obtained by solving probabilistic PF as described in Section III. In this paper, the different uncertainties considered and their probabilistic models are described as follows:

*1) Probabilistic modelling of Wind Turbine Generator:*

Power output of a wind turbine generator (WTG) varies non-linearly as per wind speed variations. However, due to intermittency of wind flow, the output of a WTG also becomes intermittent. As an accurate forecast of wind speeds at a certain interval is extremely difficult to obtain, probabilistic rendering of wind speed pattern is better suited for obtaining its estimated value and hence the estimated WTG output. Wind speed variations are widely considered to follow Weibull's distribution [5], [17-20]. Probability density function (pdf) of Weibull's distribution is provided as:

$$f_U(U) = \frac{\beta}{\alpha}(U/\alpha)^{\beta-1} \exp\left(-\left(\frac{U}{\alpha}\right)^\beta\right) \quad (16)$$

Scale ($\alpha > 0$) and shape ($\beta > 0$) parameters of the Weibull's pdf can be obtained from the mean wind speed and standard deviation. WTG output is defined as a function of input wind speed as per the following equation:

$$P_{WTG} = \begin{cases} 0 & 0 \leq U < U_{ci} \\ P_{WTG_{rt}} \times \frac{U - U_{ci}}{U_{rt} - U_{ci}} & U_{ci} \leq U < U_{rt} \\ P_{WTG_{rt}} & U_{rt} \leq U \leq U_{co} \\ 0 & U_{co} \leq U \end{cases} \quad (17)$$

2) *Probabilistic modelling of system loads:*

Loads at system buses are considered to follow normal distribution with their mean same as the base system data [5], [19]. Standard deviations of the loads ($\sigma_L > 0$) are set at a specific percentage of their mean values ($\mu_L > 0$). Its pdf is:

$$f_L(L) = \frac{1}{\sqrt{2\pi\sigma_L^2}} \exp\left(-\frac{(L-\mu_L)^2}{2\sigma_L^2}\right) \quad (18)$$

3) *Probabilistic modelling of line contingencies:*

Line contingencies are modelled in this work according to binomial distribution [19] with the line forced outage rates (FORs = $\rho_{FOR}$) providing the probability of their failures. A single trial and output values of either 1 (line in service) or 0 (line outage) is considered for modelling. It essentially becomes a Bernoulli distribution, with 1 Bernoulli trial. $H = 0, 1$.

$$f_H(H) = \binom{1}{H}(1-\rho_{FOR})^H (\rho_{FOR})^{1-H} \quad (19)$$

III. POINT ESTIMATE METHOD

Formulated by Hong, in 1998 [23], point estimate method (PEM) is used to estimate the expected value of a function of random variables from a limited number of evaluations of the function at a few points. Two PEM and three PEM are its most common variations, with the latter providing a greater accuracy than the former. It considerably reduces the computational burden compared to MCS method, where, all realizations of a random variable are required to obtain the expected value of the output function. Three PEM is used in this chapter for the development of the proposed solution methodology for PrACTNEP. A brief overview is provided here.

Let us consider a discrete single random variable $x$ (input variable), with $Y$ distinct points. So, its pdf, $f_x(x)$ becomes same as its probability distribution function, and $Z = h(x)$ (output variable) be a function of $x$. Let, the mean or expectation, standard deviation (SD) and coefficient of variation of $x$ be denoted by $\mu_x$, $\sigma_x$ and $\nu_x$ respectively. Therefore,

$$\mu_x = E(x) = \sum x f_x(x) \quad (20)$$

$$\sigma_x^2 = E(x^2) - E(x)^2 = \sum (x - \mu_x)^2 f_x(x) \quad (21)$$

$$\nu_x = \frac{\mu_x}{\sigma_x} \quad (22)$$

Similarly, mean and SD of $Z$ can be obtained by replacing $x$ by $h(x)$. Here, it is intended to obtain the actual mean of the output variable from the three points of realization of the output function $h(x)$. Therefore, we must have:

$$p_1 h(x_1) + p_2 h(x_2) + p_3 h(x_3) = E(h(x)) \quad (23)$$

Where, $p_1$, $p_2$ and $p_3$ are the weightages for the random variable locations $x_1$, $x_2$ and $x_3$ respectively. These weightages are required to be determined from three PEM, the process of which is elaborated in the Appendix.

*A. Solution of Probabilistic PF by Three PEM [24]:*

Distinct advantage of the point estimate method is the ability for its generalization. For a function $Z = h(\mathbf{G})$, $\mathbf{G} = [g^1 \, g^2 \, g^3 \, \ldots \, g^m]$, dependent on $m$ distinct random variables, $(g^1, g^2, g^3, \ldots, g^m)$, the procedure to obtain the mean $\mu_z = E(h(\mathbf{G}))$ is as follows:

1. Evaluate the mean and variances of all random variables $g^b$; $b = 1, 2, \ldots, m$.
2. Calculate $p_a^b$, $\xi_a^b$, $\lambda_{g^b,3}$, and $\lambda_{g^b,4}$ for $g^b$, $\forall b$; $a = 1, 2$. $p_3^b = \left(\frac{1}{m}\right) - p_1^b - p_2^b$. $\xi_3^b = 0$.
3. Set $E(Z) = 0$
4. Choose a random variable $g^b$, $b = 1, 2, \ldots, m$.
5. Calculate the point estimates of $g^b$.
6. Solve deterministic power flow (PF): Evaluate $Z = R(\mathbf{G})$, with $g^b$ while keeping other random variables at their mean values i.e.,
$Z = R(g_{\mu_{g^1}}^1, g_{\mu_{g^2}}^2, g_{\mu_{g^3}}^3, \ldots, g_a^b, \ldots, g_{\mu_{g^m}}^m), a = 1, 2, 3$
7. Update $E(Z) = E(Z) + p_a^b \times R(g_{\mu_{g^1}}^1, g_{\mu_{g^2}}^2, g_{\mu_{g^3}}^3, \ldots, g_a^b, \ldots, g_{\mu_{g^m}}^m)$.
8. Terminate if all random variables have been taken into consideration, else repeat from step 4.

It can be observed that, to reach a fair estimate of $\mu_z$, the deterministic solution of $Z$ is required to be solved for $3m$ number of times. Thus, for $m$ random variables, this method is known as $3m$ point estimate method. Now, due to $\xi_3^b = 0$, $m$ number of points lie at the mean, $\mathbf{G} = [g_{\mu_{g^1}}^1 \, g_{\mu_{g^2}}^2 \, \ldots, g_{\mu_{g^m}}^m]$. It is thus sufficient to solve $Z$ only once at this point, with an augmented value of weightage,

$$p_0 = \sum_{b=1}^{m} p_3^b \quad (24)$$

The $3m$ point estimate method involves $2m + 1$ repetitions of evaluating the function $Z$ to obtain $\mu_z$. Therefore, it is also known as $2m + 1$ point estimate method.

Solution of probabilistic ACPF can be obtained conveniently if $Z$ is replaced by deterministic ACPF.

IV. PROPOSED SOLUTION METHODOLOGY

Solution of PrACTNEP involves tremendous computational burden due to the requirement of solving probabilistic ACPF for each solution string generated by the metaheuristic algorithm. This contributes the most to the overall computational requirement. Therefore, straight forward solution of PrACTNEP for even a medium sized system becomes a formidable task. In this study, the primary objective is to develop suitable intelligent strategies for efficient solution of the same.

Here, a two-stage solution methodology is proposed, which produces acceptably good sub-optimal solutions with extremely low computational burdens. The proposed methodology produces fast and efficient planning results for various systems as will be evident in Section V. The key for reducing the computational burden in PrACTNEP is to reduce the number of times probabilistic PF is solved. The stages and strategies for the methodology are discussed in detail as follows:

*A. Stage 1: Use the deterministic ACTNEP solution as an initial starting point:*

Starting from a good initial guess helps in faster arrival at the final solution by MABC, thereby reducing the computational requirement. Use of the deterministic/crisp ACTNEP solution as an initial starting point is logical as it provides a fair sub-optimal solution to start with. Additionally, as $v_{pr}$ is invariably higher than $v_{cr}$, due to uncertainties, the latter provides a lower limit to the final investment cost. It is quite justified to assume that; the final cost will never be lesser than crisp TNEP cost and the combinations which suggest a lower cost can be readily removed from the solution process.

*B. Stage 2: Intelligent strategies to reduce the number of probabilistic AC power flow solutions:*

In this stage, several intelligent and logical strategies are applied to reduce the required number of times of probabilistic ACPF solutions. These are as follows:

*1) Place a bound on the number of power corridors and total number of new lines:*

Experience of solving PrACTNEP without the use of any intelligent strategies has shown that the final solution contains new lines in almost 90% of the power corridors present in AC crisp solution. It is obvious that in the solution with uncertainties, there will be more power corridors which contain new line additions. However, for the sake of generalization with different test systems, a bound of 90-130% of the number of power corridors present in deterministic ACTNEP is placed for solving PrACTNEP. Only when a combination produced by MABC has the number of new power corridors within this bound, probabilistic ACPF is solved. For other cases, a suitably high penalty is added to the combination cost so that it is discarded by MABC algorithm. This helps in filtering out numerous combinations and narrows down the search space. Although, it should be noted that if the bound is set too tight, MABC algorithm may get trapped in a local optimum. The bound used in this study has been obtained through numerous trials and it provides a proper balance of directed search and flexibility to the MABC for obtaining a good final solution.

In a similar manner, a bound on the total number of new lines is also placed. Probabilistic ACPF for a combination is solved only when the total number of its new lines lies between 70-200% of the number in the crisp solution. For the combinations where this condition is not satisfied, a suitable penalty is added similar to the previous strategy. Here also, the bound is carefully selected through several trials to maintain proper balance of algorithm search properties.

*2) Solve probabilistic AC power flow solutions only for worthy combinations:*

As with any metaheuristic, most of the initial combinations produced by MABC are infeasible. Gradually, MABC improves upon these to find feasible sets of combinations.

The only plausible way to know about these infeasibilities is to obtain their fitness values by solving probabilistic AC power flow for each combination produced. This is an extremely inefficient process and involves huge computational burden. An intelligent way to improve the efficiency is to predict the quality or fitness of a combination without actually solving power flow.

It has been found by several trials with different systems that the final PrACTNEP cost remains within twice of the crisp ACTNEP cost, $v_{cr}$. This criterion is used to determine the quality of a combination, and set an upper limit ($v_{ulim}$). Additionally, the cost of a combination cannot be less than $v_{cr}$. Whenever a combination produced has its investment cost within the bound specified by $v_{cr}$ and $v_{ulim}$, probabilistic ACPF is solved, else, a suitable penalty is added to the combination.

The upper limit ($v_{ulim}$) is set dynamically over the course of the solution process to take advantage of a feasible lower cost combination, if obtained. To start with, a fairly relaxed limit of $v_{ulim} = 2 \times v_{cr}$ is set to allow a good variance in the population pool of the algorithm. As stated earlier, here also, if a very tight limit is placed, it can result in stagnation of the solution process.

*3) Continue solving AC power flows only if feasible solutions are obtained with previous uncertain variables:*

Solving, a probabilistic ACPF requires repeated solutions of deterministic ACPF as discussed in Section III. The final aim of solving PrACTNEP is to obtain feasibility (without any limit violations) for all probabilistic realizations of the uncertain input variables. Therefore, the overall time of solving such probabilistic ACPF can be drastically reduced if the solution process is truncated once an infeasibility is obtained while solving deterministic ACPF for an uncertain variable, for any location determined by the $2m + 1$ point estimate method. After such truncation, the combination is amended with a suitable proportionate penalty depending on the uncertain variable and its location.

*4) Continue checking feasibilities for next line contingency, only if base case and previous line contingencies produce feasible solutions:*

This is a logical extension of the previous strategy. While solving full N-1 contingency constrained PrACTNEP, in addition to the strategy 3, the solution process is truncated as soon as infeasibilities are obtained in the base case or in any line contingency. Then, proportionate penalty is added similar to the previous case.

*5) RPP of the system is fixed to that of the crisp planning:*

As stated earlier, solving RPP with ACTNEP involves substantial increase in computational burden. This increases even more for PrACTNEP. Therefore, here, to keep the computational burden within manageable proportions, RPP is considered same as that of the crisp ACTNEP [27]. Such consideration is quite reasonable also, since it has been observed that, accurate RPP does not result in much reduction of line cost compared to an approximate RPP, especially when generous limits on voltage magnitude variations are considered at the system buses.

Application of these intelligent strategies along with the two-stage solution process results in drastic reduction in computational burden compared to rigorous single stage solution process, which does not use any such intelligent strategies. This will be evident from the results in the next section. A detailed flow chart of the solution methodology is provided in Fig. 1.

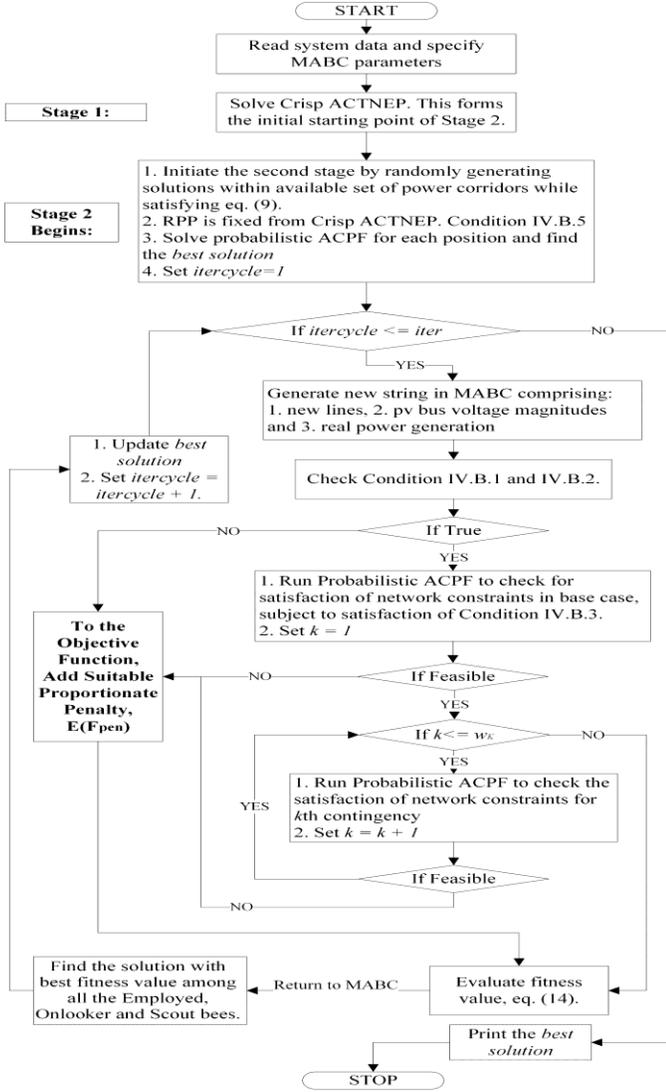

Fig. 1. Flow Chart of the Proposed Solution Methodology

## V. RESULTS AND DISCUSSION

To demonstrate the potential of the proposed methodology, it is used to solve probabilistic TNEP for Garver 6 bus [5], [28] and IEEE 24 bus [5], [28], [29] systems.

Planning problems for both the systems are solved considering dispatchable thermal generations for non-renewable resources. Wind generations present are considered to be operating under a wind pdf of Weibull distribution, with a given mean and SD. A FOR of 0.99 is considered for accurate probabilistic modelling of all possible line outages in a network. In addition, deterministic N-1 contingency in the presence all other uncertainties are also considered to obtain a full N-1 security constrained PrACTNEP. For both types of studies, system bus voltage magnitudes are limited to $\pm 5\%$ of their nominal values when line outages are not considered. In presence of line outages, a more relaxed $\pm 10\%$ limitation on bus voltage magnitudes are placed. To account for the generation and load uncertainties, conventional generators are assumed to dispatch according to their respective participation factors.

Simulation studies are conducted with MATLAB R2015b on a desktop computer having 16 GB of RAM and Intel (R) Core (TM) i5-4590 CPU processor @ 3.30 GHz. A sampling size of 10,000 is considered to transform continuous random variables to discrete. The results presented in this work are the best among 50 trials conducted for each case. Solution times described correspond to the total time required for a single trial. The settings of MABC, which are obtained by the same procedure as in [25]–[27] and used in this study, are as follows: $cs_N = 20, \psi = 2, lim = 6, iter = 30$, and $w_g = 1.5$.

### A. Probabilistic DC and ACTNEP for Garver 6 Bus System

Electrical and system data for this system is considered as in [5], for DCTNEP and from [28] for ACTNEP. Here, at bus 3, 120 MW generation is replaced by wind generation, whereas, a 10% variation of system real and reactive loads over their mean values are considered.

#### 1) Probabilistic DCTNEP:

PrDCTNEP results obtained by the proposed methodology for this system with line FOR are provided in Table 1. Here, probabilistic AC power flows for all the intelligent strategies developed in the previous section are replaced by probabilistic DC power flows. Also, the strategy IV.B.5 is not used as RPP is not required in PrDCTNEP. When compared with the results for no load shedding in [5], it can be observed that the proposed methodology is able to obtain a planning with 15.38% reduction in total cost within just a few seconds.

TABLE I
Probabilistic DCTNEP results obtained for Garver 6 bus system

|  | Proposed Method | Reported Results [5] |
|---|---|---|
| New lines Constructed | $n_{2-6} = 3; n_{3-5} = 2$; $n_{4-6} = 3$ | $n_{1-5} = 1; n_{2-3} = 1; n_{2-6} = 3$; $n_{3-5} = 2; n_{4-6} = 3$ |
| No. of New Lines | 8 | 10 |
| $v_{pr}$ (x $10^3$ US$) | **220** | 260 |
| $tp$ | 6.68 secs | - |

#### 2) Probabilistic ACTNEP:

As appropriate reference works are unavailable, solution of PrACTNEP for this system is performed by both the proposed and single-stage rigorous methods for comparison. Considering full N-1 contingency, the results are shown in Table II. It can be observed that, the application of the proposed method has resulted in a similar line planning compared to the rigorous method. However, the reduction in computational burden obtained by the former is above 99% over the latter.

### B. Probabilistic DC and ACTNEP for IEEE 24 Bus System

In this modified IEEE 24 bus system, half of the conventional generation capacities at buses 7 and 22 are replaced by wind generation, resulting in 450 MW of wind generation at each bus. Wind pdf follows Weibull distribution as utilized in [5] for obtaining PrDCTNEP. System load variations at each bus are 5% of their mean values.

TABLE II
N-1 Contingency Constrained PrACTNEP results for Garver 6 bus system

|  | Crisp Results | Proposed Method | Rigorous Method |
|---|---|---|---|
| New lines Constructed | $n_{2-6}=2$; $n_{3-5}=2$ $n_{4-6}=2$ | $n_{2-3}=2$; $n_{2-6}=3$ $n_{3-5}=2$; $n_{4-6}=3$ | $n_{2-3}=2$; $n_{2-6}=3$ $n_{3-5}=2$; $n_{4-6}=3$ |
| No. of New Lines | 6 | 10 | 10 |
| $v_{pr}$ (x $10^3$ US$) | **160.0000** | **260.0000** | **260.0000** |
| $tp$ | 21.51 secs | 104.93 secs | 3.57 hrs |
| % Reduction in time of solution when Proposed Method is used | | | **99.18** |

*1) Probabilistic DCTNEP:*

For the 24-bus network used in [5], PrDCTNEP with line FOR is solved by the proposed method and results are provided in Table III. When compared with the planning reported in [5], it can be observed that, the proposed method obtains the same planning with a tremendous reduction (more than 90%) of computational burden as is evident from its time required per trial. The result obtained by the proposed method clearly indicates that, it is an extremely efficient procedure for obtaining a proper expansion plan under uncertainties, with a drastic reduction in computational burden.

TABLE III
Probabilistic DCTNEP results obtained for IEEE 24 bus system

|  | Proposed Method | Reported Results [5] |
|---|---|---|
| New lines Constructed | $n_{1-5}=1$; $n_{2-4}=1$; $n_{3-9}=1$ $n_{3-24}=1$; $n_{6-10}=2$; $n_{7-8}=2$ $n_{8-9}=2$; $n_{9-11}=1$; $n_{10-11}=1$ $n_{10-12}=1$; $n_{10-11}=1$; $n_{10-12}=1$ $n_{11-13}=1$; $n_{12-13}=1$; $n_{14-16}=1$ $n_{15-21}=1$; $n_{15-24}=1$; $n_{16-17}=1$ $n_{17-18}=1$; $n_{20-23}=1$ | $n_{1-5}=1$; $n_{2-4}=1$; $n_{3-9}=1$ $n_{3-24}=1$; $n_{6-10}=2$; $n_{7-8}=2$ $n_{8-9}=2$; $n_{9-11}=1$; $n_{10-11}=1$ $n_{10-12}=1$; $n_{10-11}=1$; $n_{10-12}=1$ $n_{11-13}=1$; $n_{12-13}=1$; $n_{14-16}=2$ $n_{15-21}=1$; $n_{15-24}=1$; $n_{16-17}=1$ $n_{17-18}=1$; $n_{20-23}=1$ |
| No. of New Lines | 22 | 22 |
| $v_{pr}$ (x $10^3$ US$) | **9020.0000** | 9020.0000 |
| $tp$ | **5.28 min** | 57 min |
| % Reduction in Computational Burden by Proposed Method | | **90.74** |

*2) Probabilistic ACTNEP:*

Due to lack of appropriate reference works for comparison, a network as described in [28], [29] is used to solve PrACTNEP for 24 bus system. However, wind and load uncertainties are considered to be same as in the DC case. Incorporating such modifications, PrACTNEP results for both the cases of considering line FORs and with full N-1 contingency are provided in Table IV. Planning results are obtained only by the proposed method as the use of rigorous method is prohibitively computation intensive. It can be observed from the results that, the investment cost and computational burden of planning with line FOR is considerably lesser than that obtained with full N-1 contingency. This is expected, as in the latter case, all of the possible line contingencies in the system are considered without their probabilistic treatment. Such a planning provides the highest reliability for network uncertainties. However, a significantly higher investment cost and computational burden are the prices to be paid.

## VI. CONCLUSION

A two-stage solution methodology based on $2m+1$ point estimate method is proposed in this paper to efficiently solve probabilistic ACTNEP. Due to the extreme computational burden involved in solving such a problem with conventional methods, several intelligent strategies have been developed. These strategies lead to significant reduction in computational burden while achieving very good quality solutions. First stage of the solution process involves solving deterministic ACTNEP. It provides a good starting point to the MABC algorithm. Additionally, the first stage solution also allows for the development of the intelligent strategies. In the second stage, computational burden is effectively reduced by reducing the number of times probabilistic ACPF is solved. The quality of results obtained by the proposed methodology for two different test systems prove its effectiveness. It is observed from the results that, compared to the conventional and previously reported methodologies, the proposed method is able to obtain near identical results with over 90% reduction in computational burden.

TABLE IV
PrACTNEP results obtained by the proposed method for IEEE 24 bus system

|  | Crisp Results | With Line FOR | N-1 Contingency |
|---|---|---|---|
| New lines Constructed | $n_{1-5}=1$; $n_{3-9}=1$ $n_{4-9}=1$; $n_{6-10}=2$ $n_{7-8}=3$; $n_{10-11}=1$ $n_{11-13}=1$; $n_{14-16}=1$ $n_{14-23}=1$; $n_{20-23}=1$ | $n_{1-5}=1$; $n_{3-24}=1$ $n_{6-10}=2$; $n_{7-8}=2$ $n_{8-10}=1$; $n_{9-12}=1$ $n_{10-11}=2$; $n_{11-13}=1$ $n_{14-16}=1$; $n_{15-24}=1$ $n_{16-17}=1$; $n_{20-23}=1$ | $n_{1-5}=1$; $n_{2-4}=1$ $n_{3-24}=1$; $n_{6-10}=2$ $n_{7-8}=2$; $n_{8-9}=2$ $n_{10-12}=1$; $n_{12-23}=1$ $n_{14-16}=2$; $n_{15-24}=1$ $n_{16-17}=2$; $n_{17-18}=1$ $n_{20-23}=1$; $n_{21-22}=1$ |
| No. of New Lines | **13** | **15** | **19** |
| $v_{pr}$ (x $10^6$ US$) | **446** | **587** | **835** |
| $tp$ | 1189.28 secs | 2846.03 secs | 4123.69 secs |
| % Increase in time of solution when N-1 contingency is considered | | | **44.89** |

Such a drastic reduction in solving PrACTNEP is extremely beneficial as it provides a tool for solving future planning problems with more complexity, like probabilistic AC generation and transmission expansion planning which were extremely hard to solve previously.

## VII. APPENDIX: THREE-POINT ESTIMATE METHOD

Let, $M'_t(x)$ denote the $t^{th}$ order central moment of $x$,

$$M'_t(x) = \sum (x - \mu_x)^t f_x(x) \quad (25)$$

and let $\lambda_{x,t}$ denote the ratio of $M'_t(x)$ to $\sigma_x^t$,

$$\lambda_{x,t} = \frac{M'_t(x)}{\sigma_x^t} \quad (26)$$

By definition, $\lambda_{x,1} = 0$, $\lambda_{x,2} = 1$; $\lambda_{x,3}$ and $\lambda_{x,4}$, are known as coefficient of skewness and kurtosis of $x$ respectively.

For a discrete random variable, the Taylor series expansion of $h(x)$ about the mean of $x$, $\mu_x$ is given by:

$$h(x) = h(\mu_x) + \sum_{y=1}^{\infty} \frac{1}{y!} h^y(\mu_x)(x - \mu_x)^y \quad (27)$$

Here, $h^y(-)$ denotes the $y^{th}$ derivative of $h(-)$ with respect to $x$. The mean of $Z$ can be calculated as:

$$\mu_z = E(h(x)) = h(\mu_x) + \sum_{y=1}^{\infty} \left[\frac{1}{y!} h^y(\mu_x) \lambda_{x,y} \sigma_x^y\right] \quad (28)$$

Now, let, $x_a = \mu_x + \xi_a \sigma_x$, $a = 1, 2, 3$ denote the $a^{th}$ location and $\xi_a$ are the constants to be determined. Let, $p_a$ be the weightages at location $x_a$. If it is assumed that the third point is at mean, then, $\xi_3 = 0$, and we need to determine only $\xi_a$, $a = 1, 2$ and $p_a$, $a = 1, 2, 3$. Multiplying (27) by $p_a$ with $x = x_a$, $a = 1, 2, 3$ and summing them up results in:

$$p_1 h(x_1) + p_2 h(x_2) + p_3 h(x_3)$$
$$= h(\mu_x)(p_1 + p_2 + p_3)$$
$$+ \sum_{y=1}^{\infty} \left[ \frac{1}{y!} h^y(\mu_x)(p_1 \xi_1^y + p_2 \xi_2^y) \sigma_x^y \right] \quad (29)$$

As (23) holds true, therefore by comparing the first five terms of the RHS of (29) with RHS of (28), we get:

$$p_1 + p_2 + p_3 = 1; \quad p_1 \xi_1^1 + p_2 \xi_2^1 = 0; \quad p_1 \xi_1^2 + p_2 \xi_2^2 = 1$$
$$p_1 \xi_1^3 + p_2 \xi_2^3 = \lambda_{x,3}; \quad p_1 \xi_1^4 + p_2 \xi_2^4 = \lambda_{x,4} \quad (30)$$

Now, solving eq. (30), we get:

$$\xi_2 = \frac{\lambda_{x,3}}{2} + \sqrt{\left[ \lambda_{x,4} - 3 \left( \frac{\lambda_{x,3}}{2} \right)^2 \right]} \quad (31)$$

$$\xi_1 = \frac{\lambda_{x,3}}{2} - \sqrt{\left[ \lambda_{x,4} - 3 \left( \frac{\lambda_{x,3}}{2} \right)^2 \right]} \quad (32)$$

$$p_1 = \frac{1}{\xi_1(\xi_1 - \xi_2)}; \quad p_2 = \frac{-1}{\xi_2(\xi_1 - \xi_2)}; \quad p_3 = 1 - p_1 - p_2 \quad (33)$$